\documentstyle [12pt,epsf, graphicx, amssymb]{article}

\begin{document}
\def\l{\lambda}
\def\m{\mu}
\def\a{\alpha}
\def\b{\beta}
\def\g{\gamma}
\def\d{\delta}
\def\e{\epsilon}
\def\o{\omega}
\def\O{\Omega}
\def\v{\varphi}
\def\t{\theta}
\def\r{\rho}
\def\bs{$\blacksquare$}
\def\bp{\begin{proposition}}
\def\ep{\end{proposition}}
\def\bt{\begin{th}}
\def\et{\end{th}}
\def\be{\begin{equation}}
\def\ee{\end{equation}}
\def\bl{\begin{lemma}}
\def\el{\end{lemma}}
\def\bc{\begin{corollary}}
\def\ec{\end{corollary}}
\def\pr{\noindent{\bf Proof: }}
\def\note{\noindent{\bf Note. }}
\def\bd{\begin{definition}}
\def\ed{\end{definition}}
\def\C{{\mathbb C}}
\def\P{{\mathbb P}}
\def\Z{{\mathbb Z}}
\def\d{{\rm d}}
\def\deg{{\rm deg\,}}
\def\deg{{\rm deg\,}}
\def\arg{{\rm arg\,}}
\def\min{{\rm min\,}}
\def\max{{\rm max\,}}

\newcommand{\norm}[1]{\left\Vert#1\right\Vert}
\newcommand{\abs}[1]{\left\vert#1\right\vert}

\newcommand{\set}[1]{\left\{#1\right\}}
\newcommand{\setb}[2]{ \left\{#1 \ \Big| \ #2 \right\} }

\newcommand{\IP}[1]{\left<#1\right>}
\newcommand{\Bracket}[1]{\left[#1\right]}
\newcommand{\Soger}[1]{\left(#1\right)}

\newcommand{\Integer}{\mathbb{Z}}
\newcommand{\Rational}{\mathbb{Q}}
\newcommand{\Real}{\mathbb{R}}
\newcommand{\Complex}{\mathbb{C}}

\newcommand{\eps}{\varepsilon}
\newcommand{\To}{\longrightarrow}
\newcommand{\varchi}{\raisebox{2pt}{$\chi$}}

\newcommand{\E}{\mathbf{E}}
\newcommand{\Var}{\mathrm{var}}

\def\squareforqed{\hbox{\rlap{$\sqcap$}$\sqcup$}}
\def\qed{\ifmmode\squareforqed\else{\unskip\nobreak\hfil
\penalty50\hskip1em\null\nobreak\hfil\squareforqed
\parfillskip=0pt\finalhyphendemerits=0\endgraf}\fi}

\renewcommand{\th}{^{\mathrm{th}}}
\newcommand{\Dif}{\mathrm{D_{if}}}
\newcommand{\Difp}{\mathrm{D^p_{if}}}
\newcommand{\GHF}{\mathrm{G_{HF}}}
\newcommand{\GHFP}{\mathrm{G^p_{HF}}}
\newcommand{\f}{\mathrm{f}}
\newcommand{\fgh}{\mathrm{f_{gh}}}
\newcommand{\T}{\mathrm{T}}
\newcommand{\K}{^\mathrm{K}}
\newcommand{\PghK}{\mathrm{P^K_{f_{gh}}}}
\newcommand{\Dig}{\mathrm{D_{ig}}}
\newcommand{\for}{\mathrm{for}}
\newcommand{\End}{\mathrm{end}}

\newtheorem{th}{Theorem}[section]
\newtheorem{lemma}{Lemma}[section]
\newtheorem{definition}{Definition}[section]
\newtheorem{corollary}{Corollary}[section]
\newtheorem{proposition}{Proposition}[section]

\begin{titlepage}

\begin{center}

\topskip 5mm

{\LARGE{\bf {Higher derivatives of functions

\vskip 4mm

with zeros on algebraic curves}}}

\vskip 8mm

{\large {\bf G. Goldman}}

\vspace{6 mm}

{Department of Applied Mathematics, Tel Aviv University,
Tel Aviv 69978, Israel. e-mail: ggoldman@tauex.tau.ac.il}

\vspace{6 mm}

{\large {\bf Y. Yomdin}}

\vspace{6 mm}

{Department of Mathematics, The Weizmann Institute of Science,
Rehovot 76100, Israel. e-mail: yosef.yomdin@weizmann.ac.il}

\end{center}

\vspace{6 mm}
\begin{center}

{ \bf Abstract}
\end{center}

{\small Let $f: B^n \rightarrow {\mathbb R}$ be a $d+1$ times continuously differentiable function on the unit ball $B^n$, with $\max_{z\in B^n} \Vert f(z) \Vert=1$. A well-known fact is that if $f$ vanishes on a set $Z\subset B^n$ with a non-empty interior, then for each $k=1,\ldots,d+1$ the norm of the $k$-th derivative $||f^{(k)}||$ is at least $M=M(n,k)>0$.

\medskip

We show that this fact remains valid for all ``sufficiently dense'' sets $Z$ (including finite ones). The density of $Z$ is measured via the behavior of the covering numbers of $Z$. In particular, the bound $||f^{(k)}||\ge \tilde M=\tilde M(n,k)>0$ holds for each $Z$ with the box (or Minkowski, or entropy) dimension $\dim_e(Z)$ greater than $n-\frac{1}{k}$.}

\end{titlepage}

\newpage


\section{Introduction}\label{Sec:Intro}
\setcounter{equation}{0}


\smallskip


\smallskip

In this paper we continue our study (started in \cite{Yom1}, and continued recently in \cite{Yom5}-\cite{Yom7} and \cite{Gol.Yom.10}-\cite{Gol.Yom.11}, of certain very special settings of the classical Whitney's smooth extension problem. A general setting of this problem, some recent results, and references can be found in \cite{Bru.Shv,Fef,Fef.Kla,Whi1,Whi2,Whi3}).

\medskip

Let $Z\subset B^n \subset {\mathbb R}^n$ be a closed subset of the unit ball $B^n$. In \cite{Yom5}-\cite{Yom7} we looked for $C^{d+1}$-smooth functions $f:B^n\to {\mathbb R}$, vanishing on $Z$. Such $C^{d+1}$-smooth (and even $C^\infty$) functions $f$ always exist, since any closed set $Z$ is the set of zeroes of a certain $C^\infty$-smooth function.

\medskip

We normalize the extensions $f$ requiring $||f||:=\max_{B^n}|f|=1$, and ask for the minimal possible norm of the last derivative $||f^{(d+1)}||:=\max_{B^n} ||f^{(d+1)}(x)||$, where the norm $||f^{(d+1)}(x)||$ at the point $x$ is the sum of the absolute values of the $(d+1)$-st derivatives of $f$ at $x$.

\medskip

We call this minimal possible norm of the last derivative of $f$ {\it the $d$-rigidity ${\cal RG}_d(Z)$ of $Z$}. In other words, for each normalized $C^{d+1}$-smooth function $f:B^n\to {\mathbb R},$ vanishing on $Z$, we have
$$
||f^{(d+1)}||\ge {\cal RG}_d(Z),
$$
and ${\cal RG}_d(Z)$ is the maximal number with this property.

\smallskip

Recent exciting developments in the general Whitney problem (see \cite{Bru.Shv,Fef,Fef.Kla} and references therein), provide essentially a complete answer to the general Whitney extension problem in any dimension.

\smallskip

Of course, the results of \cite{Fef,Fef.Kla} provide, in particular, an algorithmic way to estimate ${\cal RG}_d(Z)$ for any finite $Z$. However, our goal in the present paper, as well as in our previous papers, related to smooth rigidity (\cite{Yom1,Yom5,Yom6,Yom7}), is somewhat different: {\it we look for an explicit answer in terms of simple, and (at least, in principle) directly computable geometric (or topological) characteristics of $Z$.}

\smallskip

First we shortly recall some of basic properties of the $d$-rigidity ${\cal RG}_d(Z)$.

\bp\label{prop:basic}
In dimension one, for any $Z\subset [-1,1]$, ${\cal RG}_d(Z)=0,$ if the cardinality $|Z|\le d$, and ${\cal RG}_d(Z)\ge \frac{(d+1)!}{2^{d+1}}$ otherwise.
\ep
Of course, this basic fact (which can be proved in many ways, in particular, via the appropriate divided finite differences) is the natural  starting point in the study of smooth rigidity.

\medskip

An immediate consequence in higher dimensions (via the obvious scaling) is the following:

\bp\label{prop:basic1}
If there is a straight line $\ell$ which passes through $(d+1)$ points of the zero set $Z$ of $f$, and through a point $z_0\in B^n$ where $|f(z_0)|=\gamma >0$, then ${\cal RG}_d(Z)\ge \frac{\gamma (d+1)!}{2^{d+1}}$.
\ep

\medskip

As an immediate consequence we get

\bp\label{prop:basic2}
For $Z$ with a non-empty interior we always have ${\cal RG}_d(Z)\ge \frac{(d+1)!}{2^{d+1}}$, for any $d$, independently of the size and the geometry of $Z$.
\ep
\pr
Just restrict $f$ to a straight line $\ell$, passing through some point $z_0$, where $|f(z_0)|=1$, and crossing the interior of $Z$.
$\square$

\medskip

Notice that in these initial examples the rigidity (under the initial scaling, with $\max_{B^n}|f|=1$), necessarily jumps: it is either zero or at least $\frac{(d+1)!}{2^{d+1}}$. In higher dimensions ${\cal RG}_d(Z)$ may attain arbitrarily small positive values (see \cite{Yom5}).

\medskip

In \cite{Yom5,Yom6} the following question was discussed: {\it Can this last property be extended to other $Z$, beyond those with a non-empty interior? In particular, is it true for sufficiently dense finite sets $Z$?}

\medskip

A partial answer was obtained in \cite{Yom7}:

\bt\label{thm:main.intro1}(\cite{Yom7})
If the box dimension $dim_e(Z)$ is greater than $n-\frac{1}{d+1}$, then
$$
{\cal RG}_d(Z)\ge M=M(n,d)>0,
$$
where the positive constant $M$ depends only on $n$ and $d$.
\et
The box (or Minkowski, or entropy ...) dimension of $Z$ is the exponent $\eta$ in the asymptotic expression for the $\e$ - covering number $M(\e,Z)$ of $Z$ with the $\e$-balls: $M(\e,Z) \sim (\frac{1}{\e})^\eta$, as $\e \to 0$.

\smallskip

In particular, the result of Theorem \ref{thm:main.intro1} provides examples of discrete, but sufficiently dense, sets $Z$, for which ${\cal R}_d(Z)$ behaves in the same way as for sets with a non-empty interior.

\medskip

A non-asymptotic version of this result, also obtained in \cite{Yom7}, provides examples of {\it finite}, but sufficiently dense, sets $Z$ with the same property.

\subsection{The present paper}\label{pres.pap}

The present paper is, in a sense, a direct continuation of \cite{Yom7}. Indeed, the main new tool introduced in \cite{Yom7} was what we call the ``test curves'': those are smooth curves $\o$, of a possibly small ``high-order curvature'', passing through the set of zeros $Z(f)$. In fact, one of the main goals of the present paper is to accurately prove some of the combinatorial properties of the chain rule for the high order derivatives of the composition $f(\o)$, used in \cite{Yom7}.

\medskip

However, as the smooth rigidity is concerned, in the present paper {\it we go just in the opposite direction to that of \cite{Yom7}}: we obtain smooth rigidity results for certain classes of ``thin'' sets $Z$. Specifically, we consider $Z$ which are inside algebraic curves in ${\mathbb R}^n$, and which are sufficiently dense there.

\medskip

Our main results - Theorem \ref{th:First.General.Ineq1} and Corollaries \ref{cor:basic1} and \ref{cor.example} below - can be considered as generalizations of Proposition \ref{prop:basic1} from the straight lines to the curves of higher degree.

\medskip

Consider polynomial parametric curves $\o$ of degree $s$ in ${\mathbb R}^n$. The curves $\o$ are given in the coordinate form by $\o(t)=(\o_1(t),\ldots,\o_n(t))$, with $t\in [-1,1]$, and with $\o_i(t)$ being polynomials in $t$ of the degree at most $s$. Denote, as usual, by $[\eta]$ the integer part of $\eta$.

\bt\label{th:First.General.Ineq1}
Let $f: B^n \rightarrow {\mathbb R}$ and $\o$ be as above, with $\o([-1,1])\subset B^n$. Put $g(t)=f(\o(t))$. Then for each $t\in [-1,1]$ we have

\be\label{eq:main.ineq}
\sum^{d+1}_{|\alpha|=[\frac{d+1}{s}]+1} ||f^{(\alpha)}(\o(t))||\ge C(n,d,s)||g^{(d+1)}(t)||.
\ee
with the positive constant $C(n,d,s)$, which is explicitly given below.

\medskip

In particular, we have

\be\label{eq:main.ineq1}
\sum^{d+1}_{|\alpha|=[\frac{d+1}{s}]+1} ||f^{(\alpha)}||\ge C(n,d,s)||g^{(d+1)}||,
\ee
where the global norms of the derivatives in the second equation (\ref{eq:main.ineq1}) are the maxima of the point-wise norms with respect to $t$.
\et

Thus for $\o$ a curve of degree $s>1$ we cannot translate the derivatives of $f(\o)$ directly to the derivatives of $f$ of the same order, as in the case of $\o$ being the straight line. But still we can bound from below the sum of the norms of the derivatives of $f$ of orders from $[\frac{d+1}{s}]+1$ to $d+1$ in terms of the derivatives of $g=f(\o)$.

\smallskip

In what follows we do not assume the zero set $Z(f)$ to sit inside $\o$. We work with the intersection $Z(f)\cap \o([-1,1]).$

\medskip

As an immediate corollary of Theorem \ref{th:First.General.Ineq1} we obtain:

\bc\label{cor:basic1}
In the assumptions as above, if the curve $\o$ crosses the set of zeros $Z(f)$ at at least $d+1$ points, and passes through a certain point $z_0$ with $|f(z_0)|\ge \gamma > 0$, then

\be\label{eq:main.ineq2}
\sum^{d+1}_{|\alpha|=[\frac{d+1}{s}]+1} ||f^{(\alpha)}||\ge C(n,d,s)\frac{\gamma(d+1)!}{2^{d+1}}.
\ee
\ec


Another corollary is as follows:

\bc\label{cor.example}
Let $f: B^n \rightarrow {\mathbb R}$ be an infinitely differentiable function on $B^n,$ and let $\o$ be as above, with $\o([-1,1])\subset B^n$ not contained in the zero set $Z(f)$. Assume that $|Z(f)\cap \o([-1,1])|=\infty$. Then there is an infinite number of the derivatives orders $m$ for which
$$
||f^{(m)}||\ge C(d,m,f,\o)\frac{(m+1)!}{2^{m+1}},
$$
with the positive constants $C(d,m,f,\o)$, which are explicitly given below.
\ec


\medskip

Let us now describe our basic approach (which we call ``test curves''). It was introduced in \cite{Yom7}, and it is further developed in the present paper.

\smallskip

If we could find a straight line $\ell$ in ${\mathbb R}^n$, passing through the point $z_0$, where the absolute value $|f(z)|$ is equal to one, and through some $d+1$ distinct points in $Z$, we could immediately get the required lower bound ${\cal RG}_d(Z)\ge \frac{(d+1)!}{2^{d+1}}$ on the $d+1$-st derivative of $f$, via Proposition \ref{prop:basic1}.

\smallskip

However, for a generic finite or discrete set $Z$ any straight line $\ell$ meets $Z$ at one or two points, at most. Instead we replace $\ell$ by a polynomial curve $\o$ of a certain degree $s>1$, and try to mimic the calculations for $\ell$. First of all, this requires additional analysis of the high order chain-rule expressions for the composition $f(\o(t))$. Second, we have to construct polynomial curves $\o$, passing through some $d+1$ distinct points in $Z$ (or we have to assume the existence of such curves). Of course, what we present below on the chain rule, is a kind of a variation around the Faa di Bruno formula. However, we consider an accurate presentation of the combinatorics of this expression as necessary for our applications in \cite{Yom7} and in this paper below. So we present what we need of the chain rule explicitly, with proofs.

\smallskip

As for the problem of finding polynomial curves $\o$, passing through ``many points'' in $Z$, we do not touch it here. One can hope that also for polynomial curves $\o$ of a higher degree one can use a kind of ``discrete integral geometry'', applied for the straight lines in \cite{Yom7}. We can hope that this approach will produce polynomial curves $\o$, passing very close to a big number of points in $Z$. In this case the rest of the machinery, developed in \cite{Yom7} will work. This includes, in particular, ``small'' smooth deformations of $\o$, forcing it to pass exactly through the required points. The question is whether the minimal required density of $Z$ will be noticeable smaller than in the case of the straight lines, described in \cite{Yom7}. Anyway, this direction is not in the scope of the present paper.

\medskip

Thus, in this paper {\it we consider zero sets $Z$ which, by the assumptions, do contain ``many points'' on algebraic curves of a prescribed degree}. Notice, that in the same way as for the straight lines, the polynomial curves $\o$ of degree $s$ generically meet $Z$ at $s$ points or less, and hence {\it the sets $Z$ considered below are typically non-generic (besides the sets $Z$ with a non-empty interior)}. Still, we believe that understanding the geometry of such non-generic sets may be instructive.

\medskip

\subsection{Parametrization of the curves $\o$}\label{Sec:Param}

Our results can be naturally extended from the polynomially parametrized curves, as above, to more general algebraic curves. Their polynomial parametrization can be provided, in particular, by the technics of the Puiseux series. As a simple example, let assume that out points sit on the cusp curve $x^2-y^3=0$. Its polynomial parametrization is given by $x=t^3, y=t^2$. Thus all the results of this paper are applicable. Beyond the Puiseux series, one can hope that also some of the results of "Smooth parametrization" may be relevant (see, e.g. [Yom], [Bourget1], [Bourget2]).

\medskip

We consider this direction as an important one, and plan to present the results separately.

\medskip
\medskip

\subsection{Organization of the paper}\label{Sec:Organization}

The paper is organized as follows: in Section \ref{Sec:Comparison.Deriv} we compare the derivatives of $f$ and of the univalent function $g=f(\o(t))$.  In particular, in Subsection \ref{Sec:Composition.Deriv} we provide the main technical ingredient of our approach - the ``combinatorial'' analysis of the high order chain-rule expressions. Here we prove all the combinatorial results on our case of the chain rule, that we need. In particular, here we provide the proofs of the results, used in \cite{Yom7}. Next, on the base of Subsection \ref{Sec:Composition.Deriv}, in Subsection \ref{Sec:the.main.ineq} we produce an inequality, relating the high-order derivatives of $f$ and $g$. Next, in Section \ref{Sec:rig.on.pol.curves} we discuss sets $Z$ possessing ``many'' points lying on polynomial curves of small degree. Finally, in this Section we combine the tools of the previous sections in order to proof the main results of the paper.



\section{Comparing derivatives of $f(x)$ and $f(\o(t))$}\label{Sec:Comparison.Deriv}
\setcounter{equation}{0}

Let $f: B^n\to R$ be a $C^{d+1}$ - smooth function. In this section we prove an inequality, which compares the high order derivatives of the function $f$, and of its restriction to a given parametrized smooth curve $\o$. More accurately, we consider $C^{d+1}$ - smooth curves $\o:[-1,1]\to B^n$, given in the coordinate form by
$$
\o(t)=(\o_1(t),\ldots, \o_n(t)),
$$
and compare the derivatives of $f$ and of the composition $f\circ \o(t)$. The ``chain rule'' expressions for higher derivatives are rather complicated, so in the next section we recall these expressions, and summarise the required facts about them.

\subsection{Derivatives of $f(\o(t))$: symbolic expressions}\label{Sec:Composition.Deriv}

We consider the composition $g(t)=f(\o(t))$ on $[-1,1]$. Thus, $g(t)$ is a $C^{d+1}$-smooth function of one real variable $t\in [-1,1]$. By the formula for the derivatives of the composition we have:

$$
g'(t)=f'(\o(t))\o'(t),
$$
$$
g''(t)=f''(\o(t))\o'(t)^2+f'(\o(t))\o''(t),
$$
$$
g'''(t)=f'''(\o(t))\o'(t)^3+3f''(\o(t))\o'(t)\o''(t)+f'(\o(t))\o'''(t).
$$
We continue this list till the fifth derivative, omitting in the further formulas the arguments $t$ and $\o(t)$:
$$
g^{(iv)}=f^{(iv)}\o'^4+6f'''\o'^2\o''+f''(3\o''^2+4\o'\o''')+f'\o^{(iv)},
$$
$$
g^{(v)}=f^{(v)}\o'^5+ 10f^{(iv)}\o'^3\o''+f'''(15\o'\o''^2+5\o'^2\o''') + f''(10\o''\o'''+5\o'\o^{(iv)})+f'\o^{(v)}.
$$
$$
..................
$$
For $n=1$ the formulas above do not need any further interpretation. However, in higher dimensions $n$, we obtain an expression which sometimes must be interpreted more accurately. This interpretation, as required by our purposes, is given by Proposition \ref{Prop:chain.rule} below.


\smallskip

As an example, consider the case $n=2, \ d=1,2,3$. Here $f=f(x,y)$ is a function of two variables, and $\o(t)$ is given by the two functions
$$
\o(t)=(x(t),y(t))=(\o_1(t),\o_2(t)).
$$
We get
$$
\frac{dg(t)}{dt}=\frac{\partial f(\o(t))}{\partial x}\o'_1+\frac{\partial f(\o(t))}{\partial y}\o'_2,
$$
$$
\frac{d^2g(t)}{dt^2}=\frac{\partial^2 f(\o(t))}{\partial x^2}(\o'_1)^2+2\frac{\partial^2 f(\o(t))}{\partial x\partial y}\o'_1\o'_2+\frac{\partial^2 f(\o(t))}{\partial y^2}(\o'_2)^2+
$$
$$
+\frac{\partial f(\o(t))}{\partial x}\o''_1+\frac{\partial f(\o(t))}{\partial y}\o''_2.
$$

\newpage

And, for the third derivative $\frac{d^3g(t)}{dt^3}$ we obtain
$$
\frac{\partial^3 f(\o(t))}{\partial x^3}(\o'_1)^3+3\frac{\partial^3 f(\o(t))}{\partial x^2\partial y}(\o'_1)^2\o'_2+3\frac{\partial^3 f(\o(t))}{\partial x\partial y^2} \o'_1(\o'_2)^2+\frac{\partial^3 f(\o(t))}{\partial y^3}(\o'_2)^3+
$$
$$
+3[\frac{\partial^2 f(\o(t))}{\partial x^2}\o'_1\o''_1+\frac{\partial^2 f(\o(t))}{\partial x\partial y}(\o''_1\o'_2+\o'_1\o''_2)+\frac{\partial^2 f(\o(t))}{\partial y^2}\o'_2\o''_2]+
$$
$$
+\frac{\partial f(\o(t))}{\partial x}\o'''_1+\frac{\partial f(\o(t))}{\partial y}\o'''_2.
$$

\medskip

Clearly, the examples above, as well as the Faa di Bruno formula, suggest that a general expression for the derivative $\frac{d^{d+1}g(t)}{dt^{d+1}}$ is a sum of the partial derivatives of $f$, multiplied by certain polynomials in the derivatives of $\o$. The specific properties of these polynomials are essential for our goals, and Proposition \ref{Prop:chain.rule} below provides this necessary information. We assume below the dimension $n$ and the highest order of the derivatives $d+1$ to be fixed.

\smallskip

We need multi-index notations, and in our case this notations are not completely standard, so we recall them:

\medskip

\noindent Below $\alpha$ is a multi-index $\alpha=(\alpha_1,\cdots,\alpha_n), \ \alpha_i \in {\mathbb N}$, $|\alpha|=\alpha_1+\ldots+\alpha_n$.

\medskip

\noindent Next, for a given $\alpha, \ |\alpha| \le d+1$ put $N_\alpha=d+2-|\alpha|$. (This is the maximal order of the derivatives of $\o$ multiplied by $\frac{\partial^\alpha f}{\partial x^\alpha}$ in the expression for the derivative $\frac{d^{d+1}g(t)}{dt^{d+1}}$. See Proposition \ref{Prop:chain.rule} below). We arrange all the derivatives of $\o$ into the vector
$$
D\o_\alpha=(\o'_1,\o''_1,\ldots,\o_1^{(N_\alpha)}, \ \o'_2,\o''_2,\ldots,\o_2^{(N_\alpha)}, \ldots, \o'_n,\o''_n,\ldots,\o_n^{(N_\alpha)}).
$$

\medskip

\noindent Now consider multi-indices

$$
\beta = \beta_\alpha = (\beta^1_1,\ldots,\beta_1^{N_\alpha}, \ \beta_2^1,\ldots,\beta_2^{N_\alpha}, \ldots, \beta_n^1,\ldots,\beta_n^{N_\alpha}), \ \ \beta_i^j=0,1,2,\ldots,
$$
which are the possible multi-indices of the powers of the derivatives of $\o$ in the monomials under consideration.

\medskip

According to our notations, we have

$$
(D\o_\alpha)^\beta=(\o'_1)^{\beta^1_1}(\o''_1)^{\beta_1^2}\cdots (\o_1^{(N_\alpha)})^{\beta_1^{N_\alpha}} \cdots (\o'_n)^{\beta_n^1}(\o''_n)^{\beta_n^2}\cdots (\o_n^{(N_\alpha)})^{\beta_n^{N_\alpha}}.
$$

\medskip

Finally, we denote by $\Sigma_\alpha$ the set of the multi-indices $\beta$, satisfying the following three "Compatibility conditions":

\medskip

\noindent 1. If for some $i=1,\ldots,n,$ we have $\alpha_i=0$,  then $\beta_i^j=0, \ j=1,\ldots,N_\alpha.$

\medskip

\noindent 2. For each $i=1,\ldots,n$ we have
$$
\sum_{1\le j \le N_\alpha}\beta_i^j=\alpha_i.
$$
In particular, this implies $|\beta|=|\alpha|$. This is the total degree of the monomial $(D\o_\alpha)^\beta$, if the derivatives of $\o$ are assigned with the degree $1$, independently of their order.

\medskip

\noindent 3. $\sigma(\beta):=\sum_{i=1,\ldots,n, \ j=1,\ldots,N_\alpha} \ j\cdot \beta_i^j=d+1$. This is the total degree of the monomial $(D\o_\alpha)^\beta$, if for $s=1,2,\ldots$ the derivatives of the order $s$ of $\o$ are assigned with the degree $s$.

\medskip
\medskip

To simplify the resulting expressions, we mostly put below
$$
f^{(\alpha)}(\o(t))= \frac{\partial^\alpha f(\o(t))}{\partial x^\alpha}, \ \ \o^{(k)}(t)=\frac{d^{k}g(t)}{dt^{k}},
$$
and so on.

\medskip

For a given $\alpha$ we write $\frac{d+1}{|\alpha|}=r+\gamma, \ r \in {\mathbb N}, \ 0\le \gamma < 1,$ and put $\kappa_\alpha=r$ if $\gamma=0,$ and $\kappa_\alpha=r+1$ otherwise.

\medskip

Now we are ready to give the explicit form of the ``high-order chain rule'' in our specific case. The following result is valid for $f$ as above, any smooth parametric curve $\o(t)$ (not necessarily polynomial), and $g(t)=f(\o(t))$.

\bp\label{Prop:chain.rule}
The derivative $g^{(d+1)}(t)$ is expressed as follows:
\be\label{deriv.explicit}
g^{(d+1)}(t)=\sum_{|\alpha|=1}^{d+1}f^{(\alpha)}(\o(t))\sum_{\beta \in \Sigma_\alpha} c_{\alpha,\beta} (D\o_\alpha)^\beta.
\ee
\noindent i. Here each monomial in $(D\o_\alpha)^\beta$ has the total degree $|\alpha|$, if the derivatives of $\o$ are assigned with the degree $1$, independently of their order. If the derivatives of order $s$ of $\o$ are assigned with the degree $s$, then each monomial in $(D\o_\alpha)^\beta$ has the total degree $d+1$.

\medskip

\noindent ii. All the coefficients $c(\alpha,\beta)$ in (\ref{deriv.explicit}) are nonzero integers. There sum is bounded by $(n+1)^{d+1}(d+1)^{(d+1)(n+2)}$, as well as each $c(\alpha,\beta)$, and the total number of the terms in (\ref{deriv.explicit}).

\medskip

\noindent iii. Each monomial $(D\o_\alpha)^\beta$ in (\ref{deriv.explicit}) necessarily contains the derivatives of $\o$ of the order at least $\kappa_\alpha$. In particular, for $|\alpha|\le d$ each monomial $(D\o_\alpha)^\beta$ necessarily contains the second derivatives of $\o$.
\ep
\pr

We apply induction with respect to the order $d+1$ of the derivative of $g(t)$. The explicit expressions above (up to $d=2$) provide the basis of the induction. Now assume that the result is true for $d$. We differentiate once more (i.e. pass from $d$ to $d+1$), according to the chain rule and the rule of differentiating the product. First, we have to show that (\ref{deriv.explicit}) for $d+1$ provides a correct indexation of the high-order derivatives of $\o(t)$ in the right hand side of (\ref{deriv.explicit}), i.e. that the new powers $\beta$ belong to $\Sigma_\alpha$.

\medskip

For a given $\alpha=(\alpha_1,\cdots,\alpha_n)$ we denote by $\bar \alpha_i, \ i=1,\ldots, n$ the multi-index
$$
\bar \alpha_i=(\alpha_1,\cdots,\alpha_{i-1},\alpha_i+1,\alpha_{i+1},\cdots,\alpha_n).
$$
We get

$$
g^{(d+2)}(t)= \frac{d}{dt} \left ( \sum_{|\alpha|=1}^{d+1}f^{(\alpha)}\sum_{\beta \in \Sigma_\alpha} c_{\alpha,\beta} (D\o_\alpha)^\beta \right ) =
$$
\be\label{deriv.explicit1}
 =\sum_{|\alpha|=1}^{d+1}(\sum_{i=1}^n  f^{(\alpha_i)} \o'_i)\cdot (\sum_{\beta \in \Sigma_\alpha} c_{\alpha,\beta} (D\o_\alpha)^\beta) + \sum_{|\alpha|=1}^{d+1}f^{(\alpha)} \cdot (\frac{d}{dt}\sum_{\beta \in \Sigma_\alpha} c_{\alpha,\beta} (D\o_\alpha)^\beta).
\ee
The first expression in (\ref{deriv.explicit1}) presents a sum of the partial derivatives of $f$, of the total order between $2$ and $d+2$, multiplied by the sum of the products of the derivatives of $\o$, modified, as explained below:

$$
\sum_{|\alpha|=1}^{d+1}(\sum_{i=1}^n  f^{(\bar \alpha_i)} \o'_i)\cdot (\sum_{\beta \in \Sigma_\alpha} c_{\alpha,\beta} (D\o_\alpha)^\beta)=
\sum_{|\alpha|=1}^{d+1}\sum_{i=1}^n  f^{(\bar \alpha_i)} \cdot \sum_{\beta \in \Sigma_\alpha} c_{\alpha,\beta} (D\o_\alpha)^{\bar \beta_i}.
$$
Here the multi-index $\bar \beta_i$ is obtained from $\beta$ by replacing $\beta_i^1$ with $\beta_i^1+1.$ Clearly we have
$$
|\bar \alpha_i|=|\alpha|+1, \ \sigma(\bar \beta_i)= \sigma(\beta)+1=d+2, \ \ i=1,\ldots,n.
$$
Indeed, in each of these sums we add exactly one to the previous expression, so the induction holds. Exactly in the same way we get for each $i=1,\ldots,n$ and $l=1,\ldots,n$
$$
\sum_{1\le j \le N_{\bar \alpha_i}}(\bar \beta_i)_l^j=(\bar \alpha_i)_l,
$$
this sum increasing by one for $l=i$, and remaining as it was for $l \ne i$.

\medskip

Of course, if $(\bar \alpha_i)_l$ remains zero, then also all $(\bar \beta_i)^l_j=0$ for $1\le j \le N_{\bar \alpha_i}$. We conclude, that for the first summand in expression (\ref{deriv.explicit1}) above the compatibility conditions 1-3 are satisfied, and so the induction step holds for this part.

\medskip

Next we analyze the second summand in expression (\ref{deriv.explicit1}):
$$
 \sum_{|\alpha|=1}^{d+1}f^{(\alpha)} \cdot (\frac{d}{dt}\sum_{\beta \in \Sigma_\alpha} c_{\alpha,\beta} (D\o_\alpha)^\beta).
$$
It is enough to consider the derivatives
$$
\frac{d}{dt}\sum_{\beta \in \Sigma_\alpha} c_{\alpha,\beta} (D\o_\alpha)^\beta = \sum_{\beta \in \Sigma_\alpha} c_{\alpha,\beta} \frac{d}{dt} (D\o_\alpha)^\beta.
$$
Let's recall that by our notations we have

$$
(D\o_\alpha)^\beta=(\o'_1)^{\beta^1_1}(\o''_1)^{\beta^1_2}\cdots (\o_1^{(N_\alpha)})^{\beta^1_{N_\alpha}} \cdots (\o'_n)^{\beta^n_1}(\o''_n)^{\beta^n_2}\cdots (\o_n^{(N_\alpha)})^{\beta^n_{N_\alpha}} =
$$
$$
= \prod_{i=1,\ldots,n, \ j=1,\ldots, N_\alpha} (\o_i^{(j)})^{\beta_i^j}.
$$
As we differentiate this product, we get a sum of the products, corresponding to the differentiation of each of the factors. For each the factor $(\o_i^{(j)})^{\beta_i^j}$ we get the new product, with the term $(\o_i^{(j)})^{\beta_i^j}$ replaced by $\beta_i^j(\o_i^{(j)})^{\beta_i^j-1}\o_i^{(j+1)}.$

\medskip

Accordingly, we obtain

$$
\frac{d}{dt} (D\o_\alpha)^\beta = \sum_{i,j} \beta_i^j (D\o_\alpha)^{\hat \beta_{i,j}},
$$
where the multi-index $\hat \beta_{i,j}$ is obtained from $\beta$ by replacing $\beta_i^j$ with $\beta_i^j-1,$ and replacing $\beta_i^{j+1}$ with $\beta_i^{j+1}+1.$

\medskip

It remains to notice that the Compatibility conditions 1,2,3 above are satisfied for $\hat \beta_{i,j}$, if they were satisfied (by the induction assumption) for $\beta$. This is clear for condition 1. Next, the only two terms in $\beta$ that change are $\beta_i^j$ and $\beta_i^{j+1}$. Their sum remains unchanged, and this implies condition 2 (since $\alpha$ remains the same).

\medskip

Now, the new weighted sum for these two terms is
$$
j(\beta_i^j-1)+(j+1)(\beta_i^{j+1}+1) = j\beta_i^j + (j+1)\beta_i^{j+1} -j + j +1=j\beta_i^j + (j+1)(\beta_i^{j+1}+1).
$$
We conclude that the total new weighted sum $\sigma(\hat \beta_{i,j})$ increases by $1$ to $d+2$, as required by Condition 3.

\medskip

This completes the first step of the proof of Proposition \ref{Prop:chain.rule}, showing that the ``combinatorial'' form of expression (\ref{deriv.explicit1}) is correct.

\medskip

Then the statement $(i)$ of Proposition \ref{Prop:chain.rule} is just the reformulations of the ``Compatibility conditions'' (2) and (3).

\medskip

Next, one can show that for all $\beta \in \Sigma_\alpha$ in (\ref{deriv.explicit1}) the coefficient $c(\alpha,\beta)$ is a nonzero integer. In fact, it is easy to show that the above described transformations of the monomials, as we pass from $d$ to $d+1$, are, essentially, invertible. Hence for each $\beta \in \Sigma_\alpha$ we can find its "pre-image" in the previous induction step, and hence the new $c(\alpha,\beta)$ is a nonzero integer (???).

\medskip

Next we bound from above the coefficients $c(\alpha,\beta)$. We show that the sum $\nu_d$ of the (positive) coefficients $c(\alpha,\beta)$ in expression (\ref{deriv.explicit1}) of Proposition \ref{Prop:chain.rule} does not exceed $(n+1)^{d+1}(d+1)^{(d+1)(n+2)}$.

\medskip

We use the same induction by $d$, as in the first part of the proof. It was shown there that as we differentiate once more, we obtain the sum of two expressions. The first is
$$
\sum_{|\alpha|=1}^{d+1}\sum_{i=1}^n  f^{(\bar \alpha_i)} \cdot \sum_{\beta \in \Sigma_\alpha} c_{\alpha,\beta} (D\o_\alpha)^{\bar \beta_i}.
$$
We see that the sum $\nu^1_{d+1}$ of the coefficients in this expression is exactly $n\nu_d$.

\medskip

The second expression is
\be\label{eq:coef.1}
\sum_{|\alpha|=1}^{d+1}f^{(\alpha)} \cdot (\frac{d}{dt}\sum_{\beta \in \Sigma_\alpha} c_{\alpha,\beta} (D\o_\alpha)^\beta),
\ee
and to bound the sum $\nu^2_{d+1}$ of the coefficients of the monomials in $\o$ in this expression it is enough to consider, for a given $\alpha$, the sums
\be\label{eq:coef.2}
\frac{d}{dt}\sum_{\beta \in \Sigma_\alpha} c_{\alpha,\beta} (D\o_\alpha)^\beta = \sum_{\beta \in \Sigma_\alpha} c_{\alpha,\beta} \frac{d}{dt} (D\o_\alpha)^\beta.
\ee
In the first part of the proof it was shown that
$$
\frac{d}{dt} (D\o_\alpha)^\beta = \sum_{i,j} \beta_i^j (D\o_\alpha)^{\hat \beta_{i,j}},
$$
By our definitions we have $\beta_i^j\le N_\alpha \le d+1$. The range of the summation is $nN_\alpha \le n(d+1)$. Therefore, the sum of the coefficients in the last expression is at most $n(d+1)^2$.

\medskip

Returning to (\ref{eq:coef.2}), replacing each $c_{\alpha,\beta}$ with $\nu_d$, and bounding (very roughly) the size $|\Sigma_\alpha|$ by $d^n$,  we get for the sum of the coefficients in (\ref{eq:coef.2}) the bound $d^n\cdot n(d+1)^2 \nu_d \le n(d+1)^{n+2}\nu_d$.

\medskip

Next, returning to (\ref{eq:coef.1}) we get, because of the summation on $\alpha$, another factor $d^n$. Altogether, we have for the sum of the coefficients in (\ref{eq:coef.1}) the upper bound $\nu^2_{d+1} \le n(d+1)^{2(n+2)}\nu_d$.

\medskip

Finally, $\nu_{d+1}=\nu^1_{d+1}+\nu^2_{d+1}\le n\nu_d + n(d+1)^{n+2}\nu_d\le (n+1)(d+1)^{n+2}\nu_d,$ and iterating this bound up to $(d+1)$ we obtain $\nu_{d+1} \le (n+1)^{d+1}(d+1)^{(d+1)(n+2)}.$ This completes the proof of the bound on the coefficients.

\medskip

It remains to prove the claim (iii) of Proposition \ref{Prop:chain.rule}. Write, as above, $\frac{d+1}{|\alpha|}=r+\gamma,$ with $r$ an integer, $0\le \gamma<1,$ to include both the cases of $d+1$ divisible or not divisible by $|\alpha|$. Clearly, in any case $\kappa_\alpha \ge r,$ since the sum of the orders of the derivatives of $\o$ entering each monomial is $d+1$, while the number of the factors in each monomial is $|\alpha|$. By the same reason, if $d+1$ is not divisible by $|\alpha|$, then $\kappa_\alpha \ge r+1,$ since otherwise the sum of the orders of the derivatives of $\o$ entering each monomial, would be strictly less than $d+1$. This completes the proof of Proposition \ref{Prop:chain.rule}. $\square$

\medskip

Two the extreme terms in expression (\ref{deriv.explicit1}) of Proposition \ref{Prop:chain.rule}, corresponding to $|\alpha|=1$ and $|\alpha|=d+1,$ have an especially simple form:

\be\label{eq:lowest.term}
\sum_{i=1}^n \frac{\partial f(\o(t))}{\partial x_i}\o_i^{(d+1)}(t).
\ee
and
\be\label{eq:highest.term}
\sum_{|\alpha|=d+1} f^{(\alpha)}(\o(t))(\o'(t))^{\alpha},
\ee
where, as usual, $\o'=(\o_1',\ldots,\o'_n), \ (\o')^{\alpha} = (\o_1')^{\alpha_1}\cdot \ldots \cdot (\o_n')^{\alpha_n}$.

\medskip

\subsection{The main inequality}\label{Sec:the.main.ineq}

Let's recall our setting: we consider nonzero $C^{d+1}$-smooth functions $f:B^n\to {\mathbb R}$, vanishing on $Z$, which we normalize, assuming $\max_{B^n}|f|=1$.

\smallskip

Next, we assume $\o$ to be a polynomial parametric curve of degree $s$ in ${\mathbb R}^n$. The curve $\o$ is given in the coordinate form by $\o(t)=(\o_1(t),\ldots,\o_n(t))$, with $\o_i(t)$ being polynomials in $t$ of the degree at most $s$. Assume that $\o([-1,1])\subset B^n$.

\medskip

The following theorem provides, in our specific case, a more accurate ``chain rule'' expression for the derivatives of $f(\o(t))$, than (\ref{deriv.explicit}) of Proposition \ref{Prop:chain.rule} above. It shows that for $\o$ being polynomial of degree $s$, the summation by the multi-orders $\alpha$ of the derivatives $f^{(\alpha)}$ starts only from $|\alpha| \sim \frac{d+1}{s+1}$.

\medskip

We want to give accurate (and an integer) summation bounds, so let's put

$$
\frac{d+1}{s+1} = p+\delta,
$$
with $p$ an integer, and $0\le \delta < 1$, and define $\eta_{d,s}$ to be equal to

$$
\eta_{d,s}=p+1=[\frac{d+1}{s+1}]+1,
$$
for both $\delta=0$ and $\delta>0$, where $[\phi]$ denotes, as usual, the integer part of $\phi$.

\bt\label{th:First.General.Ineq.0}
Let $f: B^n \rightarrow {\mathbb R}$ and $\o$ be as above. Then for each $t\in [-1,1]$ we have

\be\label{eq:main.ineq5}
g^{(d+1)}(t)=\sum_{|\alpha| \ = \ \eta_{d,s}}^{d+1} f^{(\alpha)}(\o(t))\sum_{\beta \in \Sigma_\alpha} c_{\alpha,\beta} (D\o_\alpha)^\beta.
\ee
\et
\pr
Put $\frac{d+1}{s+1}=p+\delta :=\zeta.$ Thus $d+1=\zeta (s+1)$. Therefore, for each multi-index $\alpha$, we have  $\frac{d+1}{|\alpha|}=\frac{\zeta}{|\alpha|}(s+1)$. Now if $|\alpha|\le \zeta,$ we get $\frac{d+1}{|\alpha|}\ge s+1$, which implies $\kappa_\alpha \ge s+1$, where $\kappa_\alpha$ was defined before Proposition \ref{Prop:chain.rule} above.

\medskip

Now we apply the statement (iii) of Proposition \ref{Prop:chain.rule}: each monomial in (\ref{deriv.explicit}) contains the derivatives of $\o$ of the order at least $\kappa_\alpha$. Therefore, for each $|\alpha|\le \zeta$ the corresponding monomial contains the derivatives of $\o$ of the orders higher than $s$, which are zero. Thus , in fact, (\ref{deriv.explicit}) is reduced to the sum

\be\label{eq:main.ineq5}
g^{(d+1)}(t)=\sum_{|\alpha| > \zeta}^{d+1} f^{(\alpha)}(\o(t))\sum_{\beta \in \Sigma_\alpha} c_{\alpha,\beta} (D\o_\alpha)^\beta.
\ee
Finally, the condition $|\alpha| > \zeta = \frac{d+1}{s+1} = p+\delta$, is equivalent, since $|\alpha|$ is an integer, to the condition $|\alpha|= p+1=\eta_{d,s}.$ This completes the proof of Theorem \ref{th:First.General.Ineq.0}. $\square$.

\medskip

Our next result repeats, essentially, Theorem \ref{th:First.General.Ineq1} of the Introduction, providing the explicit expression for the constant involved.

\bt\label{th:First.General.Ineq}
Let $f: B^n \rightarrow {\mathbb R}$ and $\o$ be as above. Then for each $t\in [-1,1]$ we have

\be\label{eq:main.ineq}
\sum^{d+1}_{|\alpha|=\eta_{d,s}} ||f^{(\alpha)}(\o(t))||\ge C(n,d,s)||g^{(d+1)}(t)||.
\ee
with the nonzero constants $C(n,d,s)=[C_1(n,d,s)]^{-1}$, where
$$
C_1(n,s,d)=s^{2(d+1)}(n+1)^{d+1}(d+1)^{(d+1)(n+2)}.
$$
\et
\pr
We use expression (\ref{eq:main.ineq5}) of Theorem \ref{th:First.General.Ineq.0} above:
$$
g^{(d+1)}(t)=\sum_{|\alpha|=\eta_{d,s}}^{d+1} f^{(\alpha)}(\o(t))\sum_{\beta \in \Sigma_\alpha} c_{\alpha,\beta} (D\o_\alpha)^\beta.
$$
Now, we have to bound all the sums $\sum_{\beta \in \Sigma_\alpha} c_{\alpha,\beta} (D\o_\alpha)^\beta, \ \ \eta_{d,s}\le \alpha \le{d+1},$ from above. Proposition \ref{Prop:chain.rule} provides, in particular, the upper bound
$$
(n+1)^{d+1}(d+1)^{(d+1)(n+2)}
$$
for the sum of the coefficients $c_{\alpha,\beta}$.

\medskip

It remains to bound from above the derivatives of $\o$. By our assumptions,
$$
\o([-1,1])\subset B^n,
$$
and hence each coordinate polynomial $\o_i, \ i=1,\ldots,n,$ is bounded in absolute value by $1$. Therefore, by (a strongly simplified and very rough version) of Markov's inequality, the derivatives of the order $l$ of $\o$ are bounded by $s^{2l}$.

\medskip

Since in each monomial $(D\o_\alpha)^\beta$ in (\ref{eq:main.ineq5}) we have
$$
\sigma(\beta):=\sum_{i=1,\ldots,n, \ j=1,\ldots,N_\alpha} \ j\cdot \beta_i^j=d+1,
$$
we obtain $|(D\o_\alpha)^\beta|\le s^{2(d+1)}$.

\medskip

Finally, combining the bounds above, we obtain

\be\label{eq:deriv.f.o.1}
|g^{(d+1)}(t)|\le \sum_{|\alpha|=\eta_{d,s}}^{d+1}||f^{(\alpha)}(\o(t))||\cdot C_1(n,d,s),
\ee
with
$$
C_1(n,s,d)=s^{2(d+1)}(n+1)^{d+1}(d+1)^{(d+1)(n+2)}.
$$
Dividing \ref{eq:deriv.f.o.1} by $C_1(n,s,d)$ we obtain (\ref{eq:main.ineq}) with $C(n,d,s)=\frac{1}{C_1(n,d,s)}$. This completes the proof of Theorem \ref{th:First.General.Ineq}. $\square$

\bigskip

\section{Smooth rigidity on polynomial curves}\label{Sec:rig.on.pol.curves}
\setcounter{equation}{0}

Now to produce ``Smooth Rigidity'' results via the bound above, we have to provide a lower bound for the $(d+1)$-st derivative of the one-dimensional restriction $g(t)=f(\o(t))$ of $f$ to $\o$. This can be done in many ways. The most promising, from the point of view of the Whitney problem would be to assume that {\ the values of $f$} are known at a finite number of points on the curve $\o$. However, in the present paper we assume that only the geometry of the zero set $Z(f)$ is known. Thus we provide the lower bounds for the derivatives of $g(t)=f(\o(t))$ via Proposition \ref{prop:basic}.

\medskip

Let $f$ and $\o$ be as above. We have to make some additional assumptions on $f$ and $\o$.

\smallskip

Let $Z\subset B^n$ be the set of zeros of $f$. Assume that $\o([-1,1])$ is not contained in $Z$, and therefore $g(t)=f(\o(t))$ is not identically zero. Denote by $m>0$ the maximum of $|g(t)|$ on $[-1,1]$.

\bt\label{th:main1}
Assume that the cardinality $|Z\cap \o([-1,1])|\ge d+1$. Then
\be\label{eq:main.ineq3}
\sum^{d+1}_{|\alpha|=\eta_{d,s}} ||f^{(\alpha)}||\ge \frac{m(d+1)!C(n,d,s)}{2^{d+1}}.
\ee
\et
\pr
By Proposition \ref{prop:basic} (with a proper normalization), at a certain point $t_0\in [-1,1]$ we have $|g^{(d+1)}(t_0)| \ge m\frac{(d+1)!}{2^{d+1}}$. Now we apply Theorem \ref{th:First.General.Ineq} at the point $t_0$. $\square$

\medskip

As it was mentioned in the introduction, this result provides a kind of a generalization of our starting example (Proposition \ref{prop:basic1}) to polynomial curves. A natural question is whether using the curves $\o$ of higher degree, one can bound from below (as for the straight lines), the sum of the derivatives $\sum_{|\alpha|=d+1} ||f^{(\alpha)}||$ of any given order $d+1$ separately?

\medskip

Now we pass to the last result of this paper, which is an extended version of Corollary \ref{cor.example} in the Introduction.
Let $f: B^n \rightarrow {\mathbb R}$ be an infinitely differentiable function on $B^n,$ and let $Z\subset B^n$ be the set of zeros of $f$. As above, $\o$ is a polynomially parametrized curve of degree $s$, with $\o([-1,1])\subset B^n$. We assume that $\o([-1,1])$ is not contained in $Z$, and therefore $g(t)=f(\o(t))$ is not identically zero, In particular, the maximum $m$ of $|g(t)|$ on $[-1,1]$ is strictly positive.

\medskip

Put $d_1=5s$, and define $d_j, \ j=2,3,...$ inductively, as the smallest $d$ for which  $\eta_{d_j,s}>d_{j-1}$. By this definition we immediately conclude that the intervals $I_j:=[\eta_{d_j,s},d_j+1]$ are disjoint. Their length is $\theta_j:=d_j+1-\eta_{d_j,s}$.

\medskip

We do not give here an accurate expression for the degrees $d_j$. Approximately, since $\eta_{d,s}\simeq \frac{d+1}{s+1}$, we have $d_j\simeq 5(s+1)^{j}.$

\bc\label{cor.example1}
Assume that $|Z\cap \o([-1,1])|=\infty$. Than for each $j=2,3,...$ we have
\be\label{eq:main.ineq4}
\sum_{|\alpha|\in I_j} ||f^{(\alpha)}||\ge \frac{m(d+1)!C(n,d_j,s)}{2^{d_j+1}}.
\ee
In particular, there is an infinite number of the growing derivatives orders $k_j \in I_j$, for which
$$
||f^{(k_j)}||\ge \frac{m(k+1)!C(n,d_j,s)}{\theta_j 2^{d_j+1}}.
$$
\ec
\pr
By the above, the intervals $I_j$ of the possible degrees $k_j$ do not overlap. For each $j$ we now pick $k_j\in I_j$ to be the index corresponding to one of the maximal summands in (\ref{eq:main.ineq3}). Since the length of the interval $I_j$ is $\theta_j$, by (\ref{eq:main.ineq3}) we conclude that this summand satisfies
\be\label{eq:main.ineq11}
\sum_{|\alpha|=k_j} ||f^{(\alpha)}(\o(t))||\ge \frac{m(k+1)!C(n,d_j,s)}{\theta_j 2^{d_j+1}}.
\ee
This complete the proof of Corollary \ref{cor.example1}. $\square$

\section{Smooth rigidity on near-polynomial curves}\label{Sec:rig.on.near.pol.curves}
\setcounter{equation}{0}

We expect that the results of Section \ref{Sec:rig.on.pol.curves} remain valid, if we replace the polynomial curves $\o$ of degree $s$ with the ``near-polynomial curves'' $\bar \o$, defined as follows:

\medskip

\bd\label{def:near.pol.curves}

A curve $\bar \o(t)$ is called a $\gamma$ - near-polynomial curve of degree s and of order $d$, if there exists a polynomial curve $\o(t)$ of degree $s$ such, that the following inequality is valid for the derivatives of $\bar \o(t)$:

\be\label{near-pol}
\max_{t\in [-1,1]} ||\bar \o^{k}(t)-\o^{k}(t)||\le \gamma, \ \ k\le d.
\ee
\ed
As it was mentioned above, we expect that the results of Section \ref{Sec:rig.on.pol.curves} remain valid, if we replace the polynomial curves $\o$ of degree $s$ with the ``near-polynomial curves'' $\bar \o$. The proof, essentially, uses only what was used in the proof of the results above: replacing the derivatives of $\o$, which are zero, with just ``small'' derivatives of $\bar \o$, together with the bounds on the coefficients in the chain-rule expression, produces the required estimates.


\medskip

\bibliographystyle{amsplain}

\end{document}